\documentclass[11pt]{article}
\usepackage{amsfonts,latexsym,rawfonts,amsmath,amssymb,amsthm,graphicx}
\numberwithin{equation}{section}
\newtheorem{theorem}{Theorem}[section]
\newtheorem{lem}[theorem]{Lemma}
\newtheorem{thm}[theorem]{Theorem}
\newtheorem{pro}[theorem]{Proposition}
\newtheorem{cor}[theorem]{Corollary}

\newtheorem{rem}[theorem]{Remark}

\def\endproof{$\hfill\Box$\\}

\def\div{{\rm div\,}}
\def\sym{{\rm sym\,}}

\title{Nonexistence of proper $p$-biharmonic maps and Liouville type theorems I: case of $p\geq2$}
\author{Yingbo Han, Yong Luo}

\date{ }
\begin{document}
\maketitle

\begin{abstract}
Let $u: (M, g)\to (N, h)$ be a map between Riemannian manifolds $(M, g)$ and $(N, h)$. The $p$-bienergy of $u$ is defined by $E_p(u)=\int_M|\tau(u)|^pd\nu_g$, where $\tau(u)$ is the tension field of $u$ and $p>1$. Critical points of $E_p(\cdot)$ are called $p$-biharmonic maps. In this paper we will prove  nonexistence result of proper $p$-biharmonic maps when $p\geq2$. In particular when $M=\mathbb{R}^m$, we get Liouville type results under proper integral conditions , which extend the related results of Baird, Fardoun and Ouakkas \cite{BFO}.

\vskip12pt

\end{abstract}

\section{Introduction}

In the past several decades harmonic map plays a central role in geometry and analysis.
Let $u:(M,g)\rightarrow(N,h)$ be a map between Riemannian manifolds $(M, g)$ and $(N, h)$. The energy of $u$ is defined by
   \begin{equation*}
     E(u)=\int_{M}\frac{|du|^{2}}{2}d\nu_{g},
   \end{equation*}
   where $d\nu_g$ is the volume element on $(M, g)$. The Euler-Lagrange equation of $E(\cdot)$ is
\begin{equation*}
  \tau(u)=\sum_{i=1}^{m}\{ \tilde{\nabla}_{e_{i}}du(e_{i})-du(\nabla_{e_{i}}e_{i})\}=0,
\end{equation*}
where $\tilde{\nabla}$ is the Levi-Civita connection on the pullback bundle $u^{-1}TN$ and $\{e_i\}$ is a local orthonormal frame field on $M$.

In 1983, Eells and Lemaire \cite{EL} (see also \cite{ES}) proposed to consider the bienergy functional
 \begin{equation*}
  E_{2}(u)=\int_{M}\frac{|\tau(u)|^{2}}{2}d\nu_{g},
\end{equation*}
where $\tau(u)$ is the tension field of $u$. Recall that $u$ is harmonic if $\tau(u)=0$. The Euler-Lagrange equation of $E_{2}(\cdot)$ is (\cite{Ji})
  \begin{equation*}
    \tau_{2}(u):=\tilde{\triangle}\tau(u)+\sum_{i=1}^{m}R^{N}(\tau(u),du(e_{i}))du(e_{i})=0,
  \end{equation*}
  where $\tilde{\triangle}:=Tr_g(\tilde{\nabla})^2$ and $R^N$ is the Riemannian curvature tensor of $(N,h)$. To further generalize the notion of harmonic maps, Han and Feng \cite{HF} considered the $p$-bienergy$(p>1)$ functional as follows:
\begin{equation*}
   E_{p}(u)=\int_{M}|\tau(u)|^{p}d\nu_{g}.
\end{equation*}
\begin{rem} In \cite{HF} Han and Feng defined a more general object called $F$-biharmonic maps. $p$-biharmonic maps are $F$-biharmonic maps with $F(t)=(2t)^\frac{p}{2}$.
\end{rem}

We define the $p$-bitension field of $u$ by (\cite{HF})
    \begin{equation}\label{bitensionfild}
\tau_{p}(u):=p\{\tilde{\triangle}( |\tau(u)| ^{p-2}\tau(u))+\sum_{i=1}^{m}\bigg(R^{N}
  \big( |\tau(u)| ^{p-2}\tau(u),du(e_{i})\big)du(e_{i})\bigg)\}.
    \end{equation}
The Euler-Lagrange equation of $E_p(\cdot)$ is $\tau_{p}(u)=0$ and a smooth map $u$ satisfying $\tau_{p}(u)=0$ is called a $p$-biharmonic map.

\section{Nonexistence result}
It is obvious that harmonic maps are $p$-biharmonic maps when $p\geq2$. We call $p$-biharmonic maps which are not harmonic \textbf{proper} $p$-biharmonic maps. It is natural to consider when $p$-biharmonic maps are harmonic maps. There are a lot of results in this direction when $p=2$ (cf. \cite{Chen}\cite{MO}\cite{Ou} for recent surveys). Han and Feng \cite{HF} proved that $p$-biharmonic maps from a compact (oriented) manifold into a manifold with nonpositive curvature must be harmonic. In noncompact case nonexistence results of proper isometric $p$-biharmonic maps were proved in \cite{CL}\cite{Han}\cite{HF}\cite{HZ}\cite{Luo1}\cite{LM} . In \cite{HZ} Han and Zhang proved the following result.
\begin{thm}[HZ]
Let $u: (M,g)\to(N,h)$ be a $p$-biharmonic map from a Riemannian manifold $(M,g)$ into a Riemannian manifold $(N,h)$ with non-positive sectional curvature and  $a\geq 0$ be a non-negative real constant.

(i) If
$$\int_M|\tau(u)|^{a+p}dv_g<\infty,$$
and the energy is finite, that is
$$\int_M|du|^2dv_g<\infty,$$
then $u$ is harmonic.

(ii) If $Vol(M,g)=\infty$ and
$$\int_M|\tau(u)|^{a+p}dv_g<\infty,$$
then $u$ is harmonic, where $p\geq2$.
\end{thm}
 The first aim of this paper is to generalize the above theorem by releasing the integral conditions.
\begin{thm}\label{main1}
Let $u: (M, g)\to (N, h)$ be a $p$-biharmonic map ($p\geq2$) from a complete Riemannian manifold $(M, g)$ into a Riemannian manifold $(N, h)$ of nonpositive sectional curvature and  $1\leq q\leq\infty$, $p-1<s$.

(i) If $|du|$ is bounded in $L^q(M)$ and $$\int_M|\tau(u)|^sdv_g<\infty,$$
then $u$ is harmonic.

(ii) If $Vol(M, g)=\infty$ and $$\int_M|\tau(u)|^sdv_g<\infty,$$
then $u$ is harmonic.
\end{thm}
When the target manifold has strictly negative sectional curvature, we have
\begin{thm}\label{main2}
Let $u: (M, g)\to (N, h)$ be a $p$-biharmonic map ($p\geq2$) from a complete Riemannian manifold $(M, g)$ into a Riemannian manifold $(N, h)$ of strictly negative sectional curvature and $$\int_M|\tau(u)|^sdv_g<\infty$$ for some $p-1<s$. Assume that there is a point $q\in M$ such that $ranku(q)\geq2$, then $u$ is a harmonic map.
\end{thm}
\begin{rem}
Here and in the following the rank of $u$ at a point $q\in M$ is defined by the dimension of the linear space $du(T_qM)$, where $T_qM$ is the tangent bundle of $M$ at $q$.
\end{rem}
\begin{rem}
When $p=2$, Theorem \ref{main1} and Theorem \ref{main2} were proved in \cite{Luo2}, which extended previous results of Luo \cite{Luo} , Maeta \cite{Ma} and Nakauchi et al. \cite{NUG}.
\end{rem}
 Because from Schoen and Yau's paper \cite{SY} we see that a harmonic map from a complete noncompact Riemannian manifold of nonnegative Ricci curvature to a Riemannian manifold of nonpositive sectional curvature  with $\int_M|du|^qdv_g<\infty(q>1)$  must be a constant map, as a corollary of Theorem \ref{main1} we have the following Liouville type result for $p$-biharmonic maps.
\begin{cor}
Let $u: (M, g)\to (N, h)$ be a $p$-biharmonic map($p\geq2$) from a complete Riemannian manifold $(M, g)$ with $Ric^M\geq 0$ into a Riemannian manifold $(N, h)$ of nonpositive sectional curvature such that
$$\int_M|\tau(u)|^s+|du|^qdv_g<\infty,$$
where $s>p-1$ and $q>1$. Then $u$ is a constant map.
\end{cor}
\begin{rem}
This Liouville type result was first proved when $p=2$ by Baird et al. in \cite{BFO}. Though they assumed $s=q=2$, it is easy to see from their proofs that their Liouville type result holds whenever $s>1$ and $q>1$.
\end{rem}
\subsection{Proof of Theorem \ref{main1}.}
First let's prove a lemma.
\begin{lem}\label{lem}
Assume that $u:(M,g)\to (N,h)$ is a $p$-biharmonic map ($p\geq2$) from a complete manifold $(M,g)$ to a nonpositively curved manifold $(N,h)$ and $$\int_M|\tau(u)|^sdv_g<\infty$$ for some $s>p-1$. Then $|\tau(u)|$ is a constant and moreover $\tilde{\nabla}\tau(u)=0$.
\end{lem}
\proof  Let $\epsilon>0$. A direct computation shows that
\begin{eqnarray}\label{ine1}
&&\Delta(|\tau(u)|^{2p-2}+\epsilon)^\frac{1}{2}\nonumber
\\&=&(|\tau(u)|^{2p-2}+\epsilon)^{-\frac{3}{2}}(\frac{1}{2}(|\tau(u)|^{2p-2}+\epsilon)\Delta|\tau(u)|^{2p-2}-\frac{1}{4}|\nabla|\tau(u)|^{2p-2}|^2).
\end{eqnarray}
Moreover,
\begin{eqnarray}\label{ine2}
\Delta|\tau(u)|^{2p-2}&=&2|\tilde{\nabla}(|\tau(u)|^{p-2}\tau(u))|^2+2h(|\tau(u)|^{p-2}\tau(u), \tilde{\Delta}(|\tau(u)|^{p-2}\tau(u)))\nonumber
\\&=&2|\tilde{\nabla}(|\tau(u)|^{p-2}\tau(u))|^2-2\sum_{i=1}^{m}(R^{N}
  (|\tau(u)| ^{p-2}\tau(u),du(e_{i}), du(e_{i}),|\tau(u)|^{p-2}\tau(u)))\nonumber
  \\&\geq&2|\tilde{\nabla}(|\tau(u)|^{p-2}\tau(u))|^2.
\end{eqnarray}
and
\begin{eqnarray}\label{ine3}
\frac{1}{4}|\nabla|\tau(u)|^{2p-2}|^2&=&h^2(|\tau(u)| ^{p-2}\tau(u),\tilde{\nabla}(|\tau(u)|^{p-2}\tau(u)))\nonumber
\\&\leq&|\tilde{\nabla}(|\tau(u)|^{p-2}\tau(u))|^2|\tau(u)|^{2p-2}.
\end{eqnarray}
From (\ref{ine1})-(\ref{ine3}) we see that $$\Delta(|\tau(u)|^{2p-2}+\epsilon)^\frac{1}{2}\geq0,$$ which by letting $\epsilon \to 0$ implies that $$\Delta|\tau(u)|^{p-1}\geq0.$$ Then if $\int_M|\tau(u)|^sdv_g<\infty$ for $s>p-1$, by Yau's \cite{Yau} classical $L^p$ $(p>1)$ Liouville type theorem we have that there exists a constat $c$ such that $|\tau(u)|=c$.

If $c=0$ then $\tilde{\nabla}\tau(u)=0$. If $c\neq 0$, then from the proof we see that $\tilde{\nabla}(|\tau(u)|^{p-2}\tau(u))=0$, i.e. $\tilde{\nabla}\tau(u)=0$. This completes the proof.\endproof

Now let us continue to prove Theorem \ref{main1}.  From the above lemma we see that $|\tau(u)|=c$ is a constant. Hence if $Vol(M)=\infty$, we must have $c=0$, which proves (ii) of Theorem \ref{main1}. To prove (i) of Theorem \ref{main1}, we distinguish two cases. If $c=0$, we are done. If $c\neq0$, we see that $Vol(M)<\infty$ and we will get a contradiction in the following. Define a l-form on $M$ by
$$\omega(X):=\langle du(X),\tau(u)\rangle,~(X\in TM).$$
Then we have
\begin{eqnarray*}
\int_M|\omega|dv_g&=&\int_M(\sum_{i=1}^m|\omega(e_i)|^2)^\frac{1}{2}dv_g
\\&\leq&\int_M|\tau(u)||du|dv_g
\\&\leq&c Vol(M)^{1-\frac{1}{q}}(\int_M|du|^qdv_g)^\frac{1}{q}
\\&<&\infty ,
\end{eqnarray*}
where if $q=\infty$ we denote $\|du\|_{L^\infty(M)}=(\int_M|du|^qdv_g)^\frac{1}{q}$.

In addition, we consider $-\delta\omega=\sum_{i=1}^m(\nabla_{e_i}\omega)(e_i)$:
\begin{eqnarray*}
-\delta\omega&=&\sum_{i=1}^m\nabla_{e_i}(\omega(e_i))-\omega(\nabla_{e_i}e_i)
\\&=&\sum_{i=1}^m\{\langle\tilde{\nabla}_{e_i}du(e_i),\tau(u)\rangle
-\langle du(\nabla_{e_i}e_i),\tau(u)\rangle\}
\\&=&\sum_{i=1}^m\langle \tilde{\nabla}_{e_i}du(e_i)-du(\nabla_{e_i}e_i),\tau(u)\rangle
\\&=&|\tau(u)|^2,
\end{eqnarray*}
where in the second equality we used $\tilde{\nabla}\tau(u)=0$. Now by Gaffney's theorem (\cite{Ga}, see the appendix for precise statement) and the above equality we have that
$$0=\int_M-\delta\omega=\int_M|\tau(u)|^2dv_g=c^2Vol(M),$$
which implies that $c=0$, a contradiction. Therefore we must have $c=0$, i.e. $u$ is a harmonic map. This completes the proof of Theorem \ref{main1}.
\endproof
\subsection{Proof of Theorem \ref{main2}}
By Lemma \ref{lem}, $|\tau(\phi)|=c$ is a constant. We only need to prove that $c=0$. Assume that $c\neq 0$, we will get a contradiction. By the $p$-biharmonic equation and the Weitzenb\"ock formula we have at $q\in M$:
\begin{eqnarray*}
0&=&-\frac{1}{2}\Delta|\tau(u)|^{2p-2}
\\&=&-\langle\tilde{\Delta}(|\tau(u)|^{p-2}\tau(u)), \tau(u)\rangle-|\tilde{\nabla}(|\tau(u)|^{p-2}\tau(u))|^2
\\&=&\sum_{i=1}^m\langle R^N(|\tau(u)|^{p-2}\tau(u), du(e_i), du(e_i), |\tau(u)|^{p-2}\tau(u)\rangle-|\tilde{\nabla}(|\tau(u)|^{p-2}\tau(u))|^2
\\&=&\sum_{i=1}^m\langle R^N(|\tau(u)|^{p-2}\tau(u), du(e_i), du(e_i), |\tau(u)|^{p-2}\tau(u)\rangle,
\end{eqnarray*}
where in the first and fourth equalities we used Lemma \ref{lem} twice. Since the sectional curvature of $N$ is strictly negative, we must have that $du(e_i)$ is parallel to $\tau(u)$ at $q\in M$ $\forall i$, i.e. rank$u(q)\leq1$, a contradiction. This completes the proof of Theorem \ref{main2}.
\endproof

\section{Stress energy tensor and a growth formula for $p$-biharmonic maps}
In the following we will derive Liouville type results for $p$-biharmonic maps ($p\geq2$) from the $m$-dimensional Euclidean space $\mathbb{R}^m$. To do this we need to use a formula for the stress energy tensor of $p$-biharmonic maps, introduced in \cite{HF}.

Let $u:(M,g)\to (N,h)$ be a smooth map between two Rienannian manifolds. The stress $p$-bienergy tensor of $u$ is defined by
\begin{eqnarray}
S_p(u)=[(1-p)|\tau(u)|^p+p \div h(|\tau(u)|^{p-2}\tau(u), du)]g-2p
\sym h(\tilde{\nabla}(|\tau(u)|^{p-2}\tau(u)), du),
\end{eqnarray}
where $\sym T(X,Y)$ denotes symmetrization of a 2-tensor, that is $\sym T(X,Y)=\frac{1}{2}(T(X,Y)+T(Y,X))$. We have
\begin{pro}[\cite{HF}, Theorem 4.3]\label{pro1}
For any smooth map $u:(M,g)\to (N,h)$
\begin{eqnarray}
(\div S_p(u))(X)=-h(\tau_p(u), du(X))-p(p-2)|\tau(u)|^{p-2}X(\frac{|\tau(u)|^2}{2}).
\end{eqnarray}
\end{pro}
\proof In Theorem 4.3 of \cite{HF}, let $F(t)=(2t)^\frac{p}{2}$. \endproof

In particular, if $u$ is a smooth $p$-biharmonic map  we have $$(\div S_p(u))(X)=-p(p-2)|\tau(u)|^{p-2}X(\frac{|\tau(u)|^2}{2}).$$

Let $T$ be a symmetric covariant 2-tensor on a Riemannian manifold $(M,g)$ and let $X$ be a vector field on $M$. Then
$$\div(T\rfloor X)=(\div T)(X)+\frac{1}{2}\langle\mathcal{L}_Xg, T\rangle,$$
where $(T\rfloor X)(Y):=T(X,Y)$, $\mathcal{L}$ is the Lie derivative operator and $$\langle\mathcal{L}_Xg, T\rangle=\langle\mathcal{L}_Xg(e_i,e_j)T(e_i,e_j),$$ where $\{e_i\}$ is an orthonormal basis. Integrating this formula over a compact domain $U$ with smooth boundary, we obtain
\begin{eqnarray}\label{integral}
\int_{\partial U}T(X,n)d\sigma=\int_U(\div T)(X)dv_g+\frac{1}{2}\int_U\langle\mathcal{L}_Xg, T\rangle dv_g,
\end{eqnarray}
where $n$ is the outward pointing unit normal and $d\sigma$ is the volume element along $\partial U$. From Proposition \ref{pro1}, taking $T=S_p(u)$ in the above formula we have the following growth formula.
\begin{thm}\label{thm1}
Let $u:V\subseteq \mathbb{R}^m\to (N,h)$ be a $p$-biharmonic map defined on a open subset $V$ of Euclidean space $\mathbb{R}^m$ with its canonical metric $g$. Let $B_r$ be a ball of radius of $r$ contained in $V$ and $S_r=\partial B_r$. Then we have
\begin{eqnarray}\label{growth formula}
&&[2-(1-\frac{1}{p})m]\int_{B_r}|\tau(u)|^pdx \nonumber
\\&=&-(m-2)\int_{S_r}h(|\tau(u)|^{p-2}\tau(u), du(\frac{\partial}{\partial_r}))d\sigma-(1-\frac{1}{p})r\int_{S_r}|\tau(u)|^pd\sigma \nonumber
\\&-&r\int_{S_r}\frac{\partial}{\partial_r}h(|\tau(u)|^{p-2}\tau(u), du(\frac{\partial}{\partial_r}))d\sigma+2r\int_{S_r}h(|\tau(u)|^{p-2}\tau(u), \tilde{\nabla}_{\frac{\partial}{\partial_r}}du(\frac{\partial}{\partial_r}))d\sigma\nonumber
\\&+&(p-2)r\int_{B_r}|\tau(u)|^{p-2}\frac{\partial}{\partial_r}(\frac{|\tau(u)|^2}{2})dx.
\end{eqnarray}
\end{thm}
\proof In (\ref{integral}) choose $X=r\frac{\partial}{\partial_r}$, $T=S_p(u)$ and $U=B_r$. Then we have
\begin{eqnarray}\label{ine4}
r\int_{S_r}S_p(u)(\frac{\partial}{\partial_r}, \frac{\partial}{\partial_r})d\sigma=\int_{B_r}\div S_p(u)(r\frac{\partial}{\partial_r})dv_g+\int_{B_r}\langle g, S_p(u)\rangle dv_g,
\end{eqnarray}
where we used $\mathcal{L}_{r\frac{\partial}{\partial_r}}g=2g$. By definition of $S_p(u)$  we see that
\begin{eqnarray}\label{ine5}
r\int_{S_r}S_p(u)(\frac{\partial}{\partial_r}, \frac{\partial}{\partial_r})d\sigma
&=&r\int_{S_r}(1-p)|\tau(u)|^p+p\div h(|\tau(u)|^{p-2}\tau(u),du)d\sigma \nonumber
\\&-&2rp\int_{S_r}h(\tilde{\nabla}_{\frac{\partial}{\partial_r}}|\tau(u)|^{p-2}\tau(u), du(\frac{\partial}{\partial_r}))d\sigma \nonumber
\\&=& r\int_{S_r}(1-p)|\tau(u)|^p+p\div h(|\tau(u)|^{p-2}\tau(u),du)d\sigma \nonumber
\\&-&2rp\int_{S_r}\nabla_{\frac{\partial}{\partial_r}}h(|\tau(u)|^{p-2}\tau(u), du(\frac{\partial}{\partial_r}))d\sigma \nonumber
\\&+&2rp\int_{S_r}h(|\tau(u)|^{p-2}\tau(u), \nabla_{\frac{\partial}{\partial_r}}du(\frac{\partial}{\partial_r}))\nonumber
\\&=&r\int_{S_r}(1-p)|\tau(u)|^pd\sigma\nonumber
\\&-&rp\int_{S_r}\nabla_{\frac{\partial}{\partial_r}}h(|\tau(u)|^{p-2}\tau(u), du(\frac{\partial}{\partial_r}))d\sigma \nonumber
\\&+&2rp\int_{S_r}h(|\tau(u)|^{p-2}\tau(u), \tilde{\nabla}_{\frac{\partial}{\partial_r}}du(\frac{\partial}{\partial_r})),
\end{eqnarray}
and by $$(\div S_p(u))(X)=-p(p-2)|\tau(u)|^{p-2}X(\frac{|\tau(u)|^2}{2})$$ we have
\begin{eqnarray}\label{ine6}
\int_{B_r}\div S_p(u)(r\frac{\partial}{\partial_r})dv_g=-p(p-2)\int_{B_r}|\tau(u)|^{p-2}r\frac{\partial}{\partial_r}(\frac{|\tau(u)|^2}{2})
dv_g.
\end{eqnarray}
In addition
\begin{eqnarray}\label{ine6'}
\int_{B_r}\langle g, S_p(u)\rangle dv_g&=&m(1-p)\int_{B_r}|\tau(u)|^pdv_g+mp\int_{B_r}\div h(|\tau(u)|^{p-2}\tau(u), du)dv_g \nonumber
\\&-&2p\int_{B_r}\sum_ih(\tilde{\nabla}_{e_i}|\tau(u)|^{p-2}\tau(u), du(e_i))dv_g \nonumber
\\&=&m(1-p)\int_{B_r}|\tau(u)|^pdv_g+2p\int_{B_r}|\tau(u)|^pdv_g \nonumber
\\&+&(m-2)p\int_{B_r}\div h(|\tau(u)|^{p-2}\tau(u), du)dv_g \nonumber
\\&=&(2p+m(1-p))\int_{B_r}|\tau(u)|^pdv_g \nonumber
\\&+&(m-2)p\int_{S_r}h(|\tau(u)|^{p-2}\tau(u), du(\frac{\partial}{\partial_r}))d\sigma.
\end{eqnarray}
From (\ref{ine4})-(\ref{ine6'}) we finish the proof of Theorem \ref{thm1}. \endproof

When $p=2$, this growth formula was stated (without a proof) in \cite{BFO}, where Baird et al. used this formula to prove several Liouville type theorems for biharmonic maps from $\mathbb{R}^m$. We will systematically extend their results in the next section to case of $p\geq2$.
\section{Liouville type theorem for $p$-biharmonic maps from $\mathbb{R}^m$}
We suppose in this section that $(M,g)$ is the $m$-dimensional Euclidean space $\mathbb{R}^m$ with its canonical metric. For harmonic maps with $m\neq2$, it is well known when $m=1$(\cite{SY}) and when $m\geq3$(\cite{GRSB}\cite{Sea}) that if their energy is finite, then they must be constant. This result was extended to biharmonic maps by Baird et al.(\cite{BFO}) when $m\neq4$. We will further prove such Liouville type results under proper integral conditions for general $p$-biharmonic maps when $p\geq2$.  It is a surprise that when $p>2$ we have Liouville type result in all dimensions (even if $m=2p$, when the energy functional $E_p$  is scaling invariant).

We will deal with separately the case of $m=1$ and $m\geq2$. In the later case we will use the growth formula (\ref{growth formula}), and the hypotheses is stronger.
\begin{thm}\label{m=1}
Let $u:(\mathbb{R},g)\to (N,h)$ be a $p$-biharmonic map ($p\geq2$) satisfying
\begin{eqnarray}
\int_{\mathbb{R}}(|\tau(u)|^p+|du|^p)dx<\infty.
\end{eqnarray}
Then $u$ is constant.
\end{thm}
\proof Since $u$ is $p$-biharmonic, hence we have $$\div(S_p(u))(X)=-h(\tau_p(u), du(X))-p(p-2)|\tau(u)|^{p-2}X(\frac{|\tau(u)|^2}{2})=-(p-2)X|\tau(u)|^p.
$$
Therefore
$$\div(S_p(u)(X,\cdot))=\frac{1}{2}\langle \mathcal{L}_Xg, S_p(u)\rangle-(p-2)X|\tau(u)|^p.$$
Taking $X=\frac{\partial}{\partial_x}$ and since $\mathcal{L}_{\frac{\partial}{\partial_x}}g=0$ we get $$\frac{\partial}{\partial_x}(S_p(u)(\frac{\partial}{\partial_x},\frac{\partial}{\partial_x})+(p-2)|\tau(u)|^p)=0.$$ Thus there exists a constant $C$ such that
\begin{eqnarray}\label{equ1}
S_p(u)(\frac{\partial}{\partial_x},\frac{\partial}{\partial_x})+(p-2)|\tau(u)|^p=C.
\end{eqnarray}
By definition of $S_p(u)$ we see that
\begin{eqnarray*}
&&S_p(u)(\frac{\partial}{\partial_x},\frac{\partial}{\partial_x})
\\&=&[(1-p)|\tau(u)|^p+p\div h(|\tau(u)|^{p-2}\tau(u), du)-2ph(\tilde{\nabla}_{\frac{\partial}{\partial_x}}(|\tau(u)|^{p-2}\tau(u)),du(\frac{\partial}{\partial_x}))]
\\&=&(1-p)|\tau(u)|^p+p\frac{\partial}{\partial_x}h(|\tau(u)|^{p-2}\tau(u)),du(\frac{\partial}{\partial_x}))
\\&-&2p\frac{\partial}{\partial_x}h(|\tau(u)|^{p-2}\tau(u),du(\frac{\partial}{\partial_x}))
+2ph(|\tau(u)|^{p-2}\tau(u),\tilde{\nabla}_{\frac{\partial}{\partial_x}}du(\frac{\partial}{\partial_x}))
\\&=&(1+p)|\tau(u)|^p-p\frac{\partial}{\partial_x}h(|\tau(u)|^{p-2}\tau(u),du(\frac{\partial}{\partial_x})),
\end{eqnarray*}
where in the last equality we used $$\tilde{\nabla}_{\frac{\partial}{\partial_x}}du(\frac{\partial}{\partial_x})=\tau(u).$$ Thus from (\ref{equ1}) we obtain
\begin{eqnarray}\label{equ2}
[1+(2p-2)]|\tau(u)|^p=p\frac{\partial}{\partial_x}h(|\tau(u)|^{p-2}\tau(u)),du(\frac{\partial}{\partial_x}))+C.
\end{eqnarray}
By Young's inequality we have
\begin{eqnarray*}
\int^{\infty}_{-\infty}|h(|\tau(u)|^{p-2}\tau(u)),du(\frac{\partial}{\partial_x}))|dx\leq C(p)\int_{\mathbb{R}}|\tau(u)|^p+|du|^pdx<\infty.
\end{eqnarray*}
Hence there exist sequences $\{R_n\}, \{R_n'\}$ such that $\lim_{n\to \infty}R_n=\infty$ and $\lim_{n\to \infty}R_n'=-\infty$ which satisfy
$$\lim_{n\to\infty}h(|\tau(u)|^{p-2}\tau(u))(R_n),du(\frac{\partial}{\partial_x})(R_n))=0,$$
and
$$\lim_{n\to\infty}h(|\tau(u)|^{p-2}\tau(u))(R_n'),du(\frac{\partial}{\partial_x})(R_n'))=0.$$
Therefore on integrating over (\ref{equ2}) from $R_n'$ to $R_n$ we get
\begin{eqnarray}
&&\int_{R_n'}^{R_n}[1+(2p-2)]|\tau(u)|^pdx \nonumber
\\&=&p[h(|\tau(u)|^{p-2}\tau(u))(R_n),du(\frac{\partial}{\partial_x})(R_n))\nonumber
\\&-&ph(|\tau(u)|^{p-2}\tau(u))(R_n'),du(\frac{\partial}{\partial_x})(R_n'))]+C(R_n-R_n').
\end{eqnarray}
Hence we have
\begin{eqnarray}\label{equ3}
&&C=\frac{1}{R_n-R_n'}\{\int_{R_n'}^{R_n}[1+(2p-2)]|\tau(u)|^pdx-p[h(|\tau(u)|^{p-2}\tau(u))(R_n),du(\frac{\partial}{\partial_x})(R_n))\nonumber
\\&+&ph(|\tau(u)|^{p-2}\tau(u))(R_n'),du(\frac{\partial}{\partial_x})(R_n'))\}.
\end{eqnarray}
Letting $n\to \infty$ in the above quality we get $C=0$ and we obtain
\begin{eqnarray}
&&\int_{R_n'}^{R_n}[1+(2p-2)]|\tau(u)|^pdx \nonumber
\\&=&p[h(|\tau(u)|^{p-2}\tau(u))(R_n),du(\frac{\partial}{\partial_x})(R_n))\nonumber
\\&-&ph(|\tau(u)|^{p-2}\tau(u))(R_n'),du(\frac{\partial}{\partial_x})(R_n'))].
\end{eqnarray}
Letting $n\to \infty$ again we get $\int_{\mathbb{R}}|\tau(u)|^pdx=0$, implying that $\tau(u)=0$, i.e. $u$ is a harmonic map.

Recall that for harmonic maps we have the following Bochner formula (\cite{EL78})
$$\frac{1}{2}\Delta|du|^2=|\tilde{\nabla}du|^2+\langle Ric^M\nabla u,\nabla u\rangle-\sum_{i,j}\langle Rm^N(du(e_i),du(e_j)du(e_i),du(e_j)\rangle,$$
where $\{e_i\}$ is a local orthonormal frame field on $M$. Hence when $M=\mathbb{R}$ we have $\frac{1}{2}\Delta|du|^2=|\tilde{\nabla}du|^2$. Therefore $$|du|\Delta|du|=|\tilde{\nabla}du|^2-|\nabla|du||^2\geq0,$$ which implies that $|du|$ is a subharmonic function on $\mathbb{R
}$. Then by Yau's $L^p$ Liouville type theorem (\cite{Yau}) for subharmonic functions we have $|du|$ is a constant which is zero by $\int_{\mathbb{R}}|du|^pdx<\infty$. Thus we have proved that $u$ is a constant map.
 \endproof

 When $m\geq 2$ we have
 \begin{thm}\label{m>1}
 Let $u:(\mathbb{R}^m,g)\to (N,h)$ be a $p$-biharmonic map satisfying
\begin{eqnarray}
\int_{\mathbb{R}^m}(|\tilde{\nabla}du|^p+|du|^p)dx<\infty,
\end{eqnarray}
where $m\geq2$ and $p>2$. Then $u$ is a harmonic map.

Moreover $u$ is a constant map if $m\geq 3$ and in addition we assume
$$\int_{\mathbb{R}^m}|du|^qdx<\infty,$$
where $2\leq q\leq m$.
 \end{thm}
 \proof To prove this theorem we will need use the growth formula (\ref{growth formula}). Notice that
 \begin{eqnarray}\label{equmain}
&&(p-2)\int_{B_r}|\tau(u)|^{p-2}\frac{\partial}{\partial_r}(\frac{|\tau(u)|^2}{2})dx \nonumber
\\&=&\frac{p-2}{p}\int_{B_r}\frac{\partial}{\partial_r}|\tau(u)|^{p}dx
\\&=&\frac{p-2}{p}\int_{S_r}|\tau(u)|^pd\sigma-\frac{(p-2)(m-1)}{p}\int_{B_r}\frac{|\tau(u)|^p}{|x|}dx. \nonumber
 \end{eqnarray}
 Equation (\ref{equmain}) is one of our main observations.

Then from the above equality and (\ref{growth formula}) we have
 \begin{eqnarray}\label{ine7}
 &&[2-(1-\frac{1}{p})m]\frac{1}{r}\int_{B_r}|\tau(u)|^pdx+\frac{(p-2)(m-1)}{p}\int_{B_r}\frac{|\tau(u)|^p}{|x|}dx \nonumber
\\&=&\frac{-(m-2)}{r}\int_{S_r}h(|\tau(u)|^{p-2}\tau(u), du(\frac{\partial}{\partial_r}))d\sigma-(1-\frac{1}{p})\int_{S_r}|\tau(u)|^pd\sigma \nonumber
\\&-&\int_{S_r}\frac{\partial}{\partial_r}h(|\tau(u)|^{p-2}\tau(u), du(\frac{\partial}{\partial_r}))d\sigma+2\int_{S_r}h(|\tau(u)|^{p-2}\tau(u), \tilde{\nabla}_{\frac{\partial}{\partial_r}}du(\frac{\partial}{\partial_r}))d\sigma\nonumber
\\&+&\frac{(p-2)}{p}\int_{S_r}|\tau(u)|^pd\sigma.
 \end{eqnarray}
 Since $$|\tau(u)|^p\leq C(m,p)|\tilde{\nabla}du|^p,$$ on applying the Young's inequality we get
 \begin{eqnarray}\label{ine8}
 \int_{\mathbb{R}^m}h(|\tau(u)|^{p-2}\tau(u), du(\frac{\partial}{\partial_r}))dx&\leq& C(p)\int_{\mathbb{R}^m}(|\tau(u)|^pdx+|du(\frac{\partial}{\partial_r})|^pdx\nonumber
 \\&\leq& C(m,p)\int_{\mathbb{R}^m}|\tilde{\nabla}du|^p+|du|^pdx<\infty.
 \end{eqnarray}
 Hence by Lemma 3.5 in \cite{BFO} there exists an increasing sequence of $(R_n)\to \infty$ and three positive constants $C_1, C_2$ and $C_3$ such that \begin{eqnarray}\label{ine9}
 \int_{S_{R_n}}|h(|\tau(u)|^{p-2}\tau(u), du(\frac{\partial}{\partial_r}))|dx\leq\frac{C_1}{R_n},
 \end{eqnarray}
 and
 \begin{eqnarray}\label{ine10}
 C_2\leq \ln\frac{R_n}{R_{n-1}}\leq C_3.
 \end{eqnarray}
 Furthermore from (\ref{ine8}) we have
 \begin{eqnarray}\label{ine10'}
 \lim_{n\to\infty}\int_{R_{n-1}}^{R_n}\int_{S_r}|h(|\tau(u)|^{p-2}\tau(u), du(\frac{\partial}{\partial_r}))|d\sigma dr=0.
 \end{eqnarray}
 Similarly,
 \begin{eqnarray*}
 &&\int_{\mathbb{R}^m}|h(|\tau(u)|^{p-2}\tau(u), \tilde{\nabla}_{\frac{\partial}{\partial_r}}du(\frac{\partial}{\partial_r}))|dx
 \\&\leq& C(p)\int_{\mathbb{R}^m}|\tau(u)|^p+| \tilde{\nabla}_{\frac{\partial}{\partial_r}}du(\frac{\partial}{\partial_r})|^pdx
 \\&\leq& C(m,p)\int_{\mathbb{R}^m}|\tilde{\nabla}du|^pdx<\infty.
 \end{eqnarray*}
 Therefore
 \begin{eqnarray}\label{ine11}
  \lim_{n\to\infty}\int_{R_{n-1}}^{R_n}\int_{S_r}|h(|\tau(u)|^{p-2}\tau(u), \tilde{\nabla}_{\frac{\partial}{\partial_r}}du(\frac{\partial}{\partial_r})|d\sigma dr=0.
 \end{eqnarray}
 Again since
 $$\int_{\mathbb{R}^m}|\tau(u)|^pdx\leq C(m,p)\int_{\mathbb{R}^m}|\tilde{\nabla}du|^pdx<\infty,$$
 we have
 \begin{eqnarray}\label{ine12}
 \lim_{n\to\infty}\int_{R_{n-1}}^{R_n}\int_{S_r}|\tau(u)|^pd\sigma dr=0.
 \end{eqnarray}
 Now integrating over (\ref{ine7}) from $R_{n-1}$ to $R_n$ we get
  \begin{eqnarray}\label{ine13}
 &&[2-(1-\frac{1}{p})m]\int_{R_{n-1}}^{R_n}\frac{1}{r}\int_{B_r}|\tau(u)|^pdxdr
 +\frac{(p-2)(m-1)}{p}\int_{R_{n-1}}^{R_n}\int_{B_r}\frac{|\tau(u)|^p}{|x|}dx dr\nonumber
\\&=&\int_{R_{n-1}}^{R_n}\frac{-(m-2)}{r}\int_{S_r}h(|\tau(u)|^{p-2}\tau(u), du(\frac{\partial}{\partial_r}))d\sigma dr-(1-\frac{1}{p})\int_{R_{n-1}}^{R_n}\int_{S_r}|\tau(u)|^pd\sigma dr\nonumber
\\&-&\int_{R_{n-1}}^{R_n}\int_{S_r}\frac{\partial}{\partial_r}h(|\tau(u)|^{p-2}\tau(u), du(\frac{\partial}{\partial_r}))d\sigma dr+2\int_{R_{n-1}}^{R_n}\int_{S_r}h(|\tau(u)|^{p-2}\tau(u),\tilde{\nabla}_{\frac{\partial}{\partial_r}}du(\frac{\partial}{\partial_r}))d\sigma dr\nonumber
\\&+&\frac{(p-2)}{p}\int_{R_{n-1}}^{R_n}\int_{S_r}|\tau(u)|^pd\sigma dr \nonumber
\\&=&\int_{R_{n-1}}^{R_n}\frac{1}{r}\int_{S_r}h(|\tau(u)|^{p-2}\tau(u), du(\frac{\partial}{\partial_r}))d\sigma dr-(1-\frac{1}{p})\int_{R_{n-1}}^{R_n}\int_{S_r}|\tau(u)|^pd\sigma dr\nonumber
\\&-&\int_{S_{R_n}}h(|\tau(u)|^{p-2}\tau(u), du(\frac{\partial}{\partial_r}))d\sigma+\int_{S_{R_{n-1}}}h(|\tau(u)|^{p-2}\tau(u), du(\frac{\partial}{\partial_r}))d\sigma \nonumber
\\&+&2\int_{R_{n-1}}^{R_n}\int_{S_r}h(|\tau(u)|^{p-2}\tau(u), \tilde{\nabla}_{\frac{\partial}{\partial_r}}du(\frac{\partial}{\partial_r}))d\sigma dr
+\frac{(p-2)}{p}\int_{R_{n-1}}^{R_n}\int_{S_r}|\tau(u)|^pd\sigma dr,
 \end{eqnarray}
 where in the second equality we used the following computations
 \begin{eqnarray}\label{equ4}
 &&\int_{R_{n-1}}^{R_n}\int_{S_r}\frac{\partial}{\partial_r}h(|\tau(u)|^{p-2}\tau(u), du(\frac{\partial}{\partial_r}))d\sigma dr\nonumber
 \\&=&\int_{R_{n-1}}^{R_n}\int_{S_r}\frac{\partial}{\partial_r}(h(|\tau(u)|^{p-2}\tau(u), du(\frac{\partial}{\partial_r}))r^{m-1})d\omega dr\nonumber
 \\&-&(m-1)\int_{R_{n-1}}^{R_n}\frac{1}{r}\int_{S_r}h(|\tau(u)|^{p-2}\tau(u), du(\frac{\partial}{\partial_r})d\sigma dr
 \\&=&\int_{R_{n-1}}^{R_n}\frac{\partial}{\partial_r}\int_{S_r}h(|\tau(u)|^{p-2}\tau(u), du(\frac{\partial}{\partial_r}))d\sigma dr\nonumber
 \\&-&(m-1)\int_{R_{n-1}}^{R_n}\frac{1}{r}\int_{S_r}h(|\tau(u)|^{p-2}\tau(u), du(\frac{\partial}{\partial_r})d\sigma dr\nonumber
 \\&=&\int_{S_{R_n}}h(|\tau(u)|^{p-2}\tau(u), du(\frac{\partial}{\partial_r}))d\sigma-\int_{S_{R_{n-1}}}h(|\tau(u)|^{p-2}\tau(u), \nonumber du(\frac{\partial}{\partial_r}))d\sigma
 \\&-&(m-1)\int_{R_{n-1}}^{R_n}\frac{1}{r}\int_{S_r}h(|\tau(u)|^{p-2}\tau(u), du(\frac{\partial}{\partial_r})d\sigma dr.\nonumber
 \end{eqnarray}
 Therefore we get
 \begin{eqnarray}\label{ine14}
&& |\frac{(p-2)(m-1)}{p}(R_n-R_{n-1})\int_{B_{R_{n-1}}}\frac{|\tau(u)|^p}{|x|}dx dr|\nonumber
\\&\leq& |[2-(1-\frac{1}{p})m]|\int_{R_{n-1}}^{R_n}\frac{1}{r}\int_{B_r}|\tau(u)|^pdxdr \nonumber
\\&+&\frac{C}{R_{n-1}}\int_{R_{n-1}}^{R_n}\int_{S_r}|h(|\tau(u)|^{p-2}\tau(u), du(\frac{\partial}{\partial_r}))|d\sigma \nonumber
\\&+&C\int_{R_{n-1}}^{R_n}\int_{S_r}|\tau(u)|^pd\sigma dr+\int_{S_{R_n}}|h(|\tau(u)|^{p-2}\tau(u), du(\frac{\partial}{\partial_r}))|d\sigma \nonumber
\\&+&\int_{S_{R_{n-1}}}|h(|\tau(u)|^{p-2}\tau(u), du(\frac{\partial}{\partial_r}))|d\sigma \nonumber
+2\int_{R_{n-1}}^{R_n}\int_{S_r}|h(|\tau(u)|^{p-2}\tau(u), \tilde{\nabla}_{\frac{\partial}{\partial_r}}du(\frac{\partial}{\partial_r}))|d\sigma dr\nonumber
\\&+&\frac{(p-2)}{p}\int_{R_{n-1}}^{R_n}\int_{S_r}|\tau(u)|^pd\sigma dr.
 \end{eqnarray}
In addition,
 \begin{eqnarray}\label{ine15}
&& \int_{R_{n-1}}^{R_n}\frac{1}{r}\int_{B_r}|\tau(u)|^pdxdr \nonumber
\\&\leq&\int_{B_{R_n}}|\tau(u)|^pdx\int_{R_{n-1}}^{R_n}\frac{1}{r}dr\nonumber
\\&=&\ln\frac{R_n}{R_{n-1}}\int_{B_{R_n}}|\tau(u)|^pdx\leq C_3\int_{B_{R_n}}|\tau(u)|^pdx.
 \end{eqnarray}
 Then by (\ref{ine9}), (\ref{ine10'}), (\ref{ine11}), (\ref{ine12}) and (\ref{ine15}), letting $n\to \infty$, we see that the right hand of (\ref{ine14}) is bounded by $C_3\int_{\mathbb{R}^m}|\tau(u)|^pdx<\infty$, but the left hand side goes to $\infty$ since $\lim_{n\to\infty}(R_n-R_{n-1})=\infty$ by $\ln\frac{R_n}{R_{n-1}}\geq C_2>0$, if $|\tau(u)|$ dose not vanish anywhere. That is we have proved that $u$ is a harmonic map.

 Furthermore if $m\geq 3$ and $\int_{\mathbb{R}^m}|du|^qdx<\infty$ $(2\leq q\leq m)$ we have $u$ is a constant map by the well known result of \cite{GRSB} and \cite{Sea} when $q=2$, of \cite{Pri} when $2\leq q<m$ and of \cite{NT} when $q=m$. \endproof

From Theorem \ref{m>1} we can obtain the following Liouville type result.
\begin{thm}
 Let $u:(\mathbb{R}^m,g)\to (N,h)$ be a $p$-biharmonic map satisfying
\begin{eqnarray}
\int_{\mathbb{R}^m}(|\tilde{\nabla}du|^p+|du|^p)dx<\infty,
\end{eqnarray}
where $m\geq2$ and $p>2$. If $ranku(x)\leq1$, $\forall x\in \mathbb{R}^m$, $u$ is a constant map.

In particular, if  $u:(\mathbb{R}^m,g)\to (N^1,h)$ is a $p$-biharmonic map satisfying
\begin{eqnarray}
\int_{\mathbb{R}^m}(|\tilde{\nabla}du|^p+|du|^p)dx<\infty,
\end{eqnarray}
where $m\geq2$ and $p>2$. Then $u$ is a constant map.
\end{thm}
\proof From Theorem \ref{m>1} we see that $u$ is a harmonic map. To prove that $u$ is a constant map we follow the argument given at the last lines of the proof of Theorem \ref{m=1}.

Since $u$ is a harmonic map, we have the following Bochner's formula
$$\frac{1}{2}\Delta|du|^2=|\tilde{\nabla}du|^2+\langle Ric^M\nabla u,\nabla u\rangle-\sum_{i,j}\langle Rm^N(du(e_i),du(e_j)du(e_i),du(e_j)\rangle,$$
where $\{e_i\}$ is a local orthonormal frame field on $M$. Hence when $M=\mathbb{R}^m$ and $ranku\leq1$ we have $\frac{1}{2}\Delta|du|^2=|\tilde{\nabla}du|^2$. Therefore $|du|\Delta|du|=|\tilde{\nabla}du|^2-|\nabla|du||^2\geq0$, which implies that $|du|$ is a subharmonic function on $\mathbb{R
}^m$. Then by Yau's $L^p$ Liouville type theorem(\cite{Yau}) for subharmonic functions we have $|du|$ is a constant which is zero by $\int_{\mathbb{R}^m}|du|^pdx<\infty$. Thus we have proved that $u$ is a constant map.
 \endproof

 \vspace{0.5cm}

\textbf{Acknowledgement:} Yingbo Han was supported by NSF of China (No.11971415) and Nanhu Scholars Program for Young Scholars of XYNU and the Universities Young Teachers Program of Henan Province (2016GGJS-096). Yong Luo was supported by the NSF of China (No.11501421). Both authors would like to thank the anonymous reviewer for the suggestions which make this paper more readable.
 \section{Appendix}
 \begin{thm}[Gaffney's theorem]
Let $(M, g)$ be a complete Riemannian manifold. If a $C^1$ 1-form $\omega$ satisfies
that
$$\int_M|\omega|dv_g+\int_M|\delta\omega| dv_g<\infty,$$
or equivalently, a $C^1$ vector field $X$ defined by $\omega(Y) = \langle X, Y \rangle, (\forall Y \in TM)$ satisfies that
$$\int_M|X|dv_g+\int_M|div X|dv_g<\infty,$$
then $$\int_M\delta\omega dv_g=\int_Mdiv Xdv_g=0.$$
\end{thm}

\section{Ethical statements}
 The authors declare that they have no conflict of interest. This article does not contain any studies with human participants or animals performed by any of the authors. Informed consent was obtained from all individual participants included in the study.

\vspace{1cm}\sc
\hspace*{4mm}{ Yingbo Han}

School of Mathematics and Statistics, Xinyang Normal University, Xinyang, 464000, Henan, China.

{\tt yingbohan@163.com}\\

Yong Luo

School of Mathematics and statistics, Wuhan University, Wuhan 430072, China.

{\tt yongluo@whu.edu.cn}
\end{document}